\documentclass[12pt, oneside, a4paper]{article}

\usepackage{amsmath}
\usepackage{amsfonts}
\usepackage{amssymb}
\usepackage{amsthm,mathrsfs}
\usepackage{enumerate}
\usepackage{graphicx}
\newtheorem{theorem}{Theorem}[section]
\newtheorem{lemma}[theorem]{Lemma}

\newtheorem{corollary}[theorem]{Corollary}

\theoremstyle{definition}
\newtheorem{remark}[theorem]{Remark}

\newtheorem{definition}[theorem]{Definition}

\title{\textbf{Littlewood-Richardson coefficients and the eigenvalues of integral line graphs}}
\author{Mahdi Ebrahimi\footnote{ m.ebrahimi.math@ipm.ir}
 \\
 {\small\em  School of Mathematics, Institute for Research in Fundamental Sciences (IPM)},\\{\small\em P.O. Box: 19395--5746, Tehran, Iran}\\
 	\\{\small\em Declarations of interest: none}
}
\date{}

\begin{document}

\maketitle


\begin{abstract}
We first describe a system of inequalities (Horn's inequalities) that characterize eigenvalues of sums of Hermitian matrices.
 When we apply this system for integral Hermitian matrices, one can directly test it by using Littlewood-Richardson coefficients.
  In this paper, we apply Horn's inequalities to analysis the eigenvalues of an integral line graph $G$ of a connected bipartite graph.
Then we show that the diameter of $G$ is at most $2\omega(G)$, where $\omega(G)$ is the clique number of $G$. Also using Horn's inequalities, we show that for every odd integer $k\geq 19$, a non-complete $k$-regular Ramanujan graph has an eigenvalue less than $-2$.
 \end{abstract}
\noindent {\bf{Keywords:}}  Littlewood-Richardson coefficient, eigenvalue, Line graph, integral graph, Ramanujan graph. \\
\noindent {\bf AMS Subject Classification Number:}  05E10, 05C50, 05C76.

\section{Introduction}
$\noindent$ Graphs considered in this paper are simple and finite. We use \cite{spect}  as a source for undefined terms and notations. Let $\Gamma$ be a graph of \textit{order} $\nu$  and \textit{size} $e$.
   The \textit{eigenvalues} of $\Gamma$ are the eigenvalues $\gamma_i$ of its adjacency matrix $A$, indexed so that $\gamma_1 \geq \gamma_2 \geq \dots\geq \gamma_\nu$. The greatest eigenvalue, $\gamma_1$, is also called the \textit{spectral radius}.
  Since all eigenvalues of $\Gamma$ can be found by computing the eigenvalues of each component of $\Gamma$, we assume throughout the paper that $\Gamma$ is connected. Thus by the Perron-Frobenius theorem, $\gamma_1>\gamma_i$ for all $i\geq 2$.
The eigenvalues of a graph are related to many of its properties and key parameters.
 The most studied eigenvalues have been the spectral radius $\gamma_1$ (in connection with the chromatic number, the independence number and the clique number of the graph \cite{Ha},\cite{Ho},\cite{Ni},\cite{Wi}), $\gamma_2$ (in connection with the expansion property of the graph \cite{Hor}) and $\gamma_\nu$ (in connection with the chromatic and the independence number of the graph \cite{Ho} and the maximum cut \cite{Kr}). Let $\mu$ be the minimal polynomial of $A$. Then the\textit{ Hoffman Polynomial} $H$ given by $H(x)=\nu \mu(x)/\mu(\gamma_1)$ characterizes regularity of $\Gamma$ by the condition $H(A)=J$, the all-ones matrix (see \cite{Hof}). We refer the reader to the monographs \cite{Ch},\cite{Cv},\cite{Cve},\cite{Go} as well as the surveys \cite{Hor},\cite{Kr} for more details about eigenvalues of graphs and their applications.

 What can be said about the eigenvalues of a sum of two Hermitian (or real symmetric) matrices, in terms  of the eigenvalues of the summands? If $A$, $B$ and $C$ are Hermitian $n$ by $n$ matrices, we denote the eigenvalues of $A$ by $\alpha:=(\alpha_1, \alpha_2,\dots, \alpha_n)$, where $\alpha_1\geq \alpha_2\geq\dots\geq \alpha_n$, and similarly write $\beta$ and $\gamma$ for the eigenvalues (spectra) of $B$ and $C$. The question becomes: What $\alpha, \beta$ and $\gamma$ can be the eigenvalues of $n$ by $n$ Hermitian (or real symmetric) matrices $A,B$ and $C$, with $C=A+B?$ There is one obvious necessary condition, that the trace of $C$ be the sum of the traces of $A$ and $B$:
 \begin{equation}
 \sum_{i=1}^n \gamma_i=\sum_{i=1}^n \alpha_i+\sum_{i=1}^n \beta_i.
 \end{equation}
 Other conditions (some inequalities) were found by Horn \cite{H2}, all having the form
 \begin{equation}
 \sum_{k\in K} \gamma_k\leq\sum_{i\in I} \alpha_i+\sum_{j\in J} \beta_j,
 \end{equation}
 for certain subsets $I,J,K$ of $\{1, \dots,n\}$ of the same cardinality $r$, with $r<n$. We write the subsets in increasing order, so $I=\{i_1<i_2<\dots<i_r\}$, $J=\{j_1<\dots<j_r\}$ and $K=\{k_1<\dots<k_r\}$.
 In \cite{H2}, Horn undertook a systematic study of such inequalities. In fact, he prescribed sets of triples $(I, J,K)$, and he conjectured that the inequalities (2) for these triples would give both necessary and sufficient conditions for a triple $(\alpha,\beta,\gamma)$ to arise as eigenvalues of Hermitian matrices $A,B$ and $C$ with $C=A+B$.

 Horn defined sets $T^n_r$ of triples $(I,J,K)$ of subsets of $\{1,\dots,n\}$ of the same cardinality $r$, by the following inductive procedure. Set
 \begin{equation}
 U^n_r:=\{(I,J,K)|\sum_{i\in I}i+\sum_{j\in J} j=\sum_{k\in K} k+r(r+1)/2\}.
 \end{equation}
 All the triples that we have listed are in $U^n_r$. When $r=1$, set $T^n_1=U^n_1$. In general,
 \begin{equation*}
 T^n_r=\{(I,J,K)\in U^n_r|\;for\;all\;p<r\;and\;all\;(F,G,H)\;in\;T^r_p,
  \end{equation*}
  \begin{equation}
 \sum_{f\in F}i_f+\sum_{g\in G} j_g\leq\sum_{h\in H} k_h+p(p+1)/2\}.
 \end{equation}
\textbf{ Conjecture} (Horn's conjecture) \textit{A triple $(\alpha, \beta, \gamma)$ occurs as eigenvalues of Hermitian $n$ by $n$ matrices $A,B$ and  $C=A+B$, if and only if $\sum_{i=1}^n \gamma_i=\sum_{i=1}^n \alpha_i+\sum_{i=1}^n \beta_i$ and the inequalities (2) hold for every $(I,J,K)$ in $T^n_r$, for all $r<n$.}\\
$\noindent$ Lidskii announced a proof of Horn's conjecture in 1982 \cite{L2}.

\begin{theorem}\label{horn}
Horn's conjecture is true.
\end{theorem}

 Now suppose $\alpha,\beta$ and $\gamma$ are integral. Applying Cayley–Hamilton theorem, each of $\alpha, \beta$ and $\gamma$ can be considered as a \textit{partition} (a weakly decreasing finite sequence of non-negative integers). Knutson and Tao \cite{Kn} showed that  one can directly test the occurrence of the triple $(\alpha, \beta, \gamma)$ as eigenvalues of Hermitian $n$ by $n$ matrices $A,B$ and $C=A+B$ by using Littlewood-Richardson coefficients. The coefficients arising in the outer tensor product (\textit{Littlewood-Richardson coefficients}) of irreducible characters of the symmetric groups (or equivalently of Schur functions) has been of central interest in representation theory and algebraic combinatorics since the landmark paper of Littlewood and Richardson \cite{Li}. More recently, these coefficients have provided the centrepiece of geometric complexity theory in an approach that seeks to settle the $P$ versus $NP$ problem \cite{Mul}; it was recently shown to require not only positivity but precise information on the coefficients \cite{Bu}. The Littlewood-Richardson rule provides an efficient positive combinatorial description for their computation (see Sect. 2 below).

\begin{theorem}\label{good}
The Littlewood-Richardson coefficient $c^\gamma_{\alpha\beta}$ is positive exactly when $\sum_{i=1}^n \gamma_i=\sum_{i=1}^n\alpha_i+\sum_{i=1}^n \beta_i$ and the inequalities (2) are valid for all $(I,J,K)$ in $T^n_r$, and all $r<n$.
\end{theorem}
The study of integral graphs was first proposed in 1973 by Harary and Schwenk \cite{Har}. Integral graphs are very rare and difficult to be found. There are comparatively huge classes of graphs containing a very restricted number of integral graphs. For example, if we regard only graphs with a given maximum vertex degree, we get that the number of such integral graphs is finite \cite{C}. The goal of this paper is to analyses the eigenvalues of integral line graphs of bipartite graphs via Horn's inequalities. Note that the line graph of a regular complete bipartite graph is integral \cite[Theorem 1.2.16]{spect}. This shows that the number of integral line graphs of bipartite graphs is not finite.

Assume that $\lambda=(\lambda_1, \dots, \lambda_l)$ is a partition. We call $|\lambda|=\lambda_1+\dots+\lambda_l$ the \textit{size} of $\lambda$, and the number of positive parts of $\lambda$, denoted by $l(\lambda)$, the \textit{length} of $\lambda$.  We also use the notation $\mathrm{k}(\lambda)$ for the number of distinct positive parts of $\lambda$.

\begin{definition}\label{def}
For a positive integer $e$, assume that  $\alpha=(\alpha_1, \dots, \alpha_m)$ and $\beta=(\beta_1, \dots, \beta_n)$ are partitions of size $e$. Also let  $\nu:=m+n$. We define  $\mathrm{P}(\alpha,\beta)$ to be the set of all partitions $\gamma$ of size $2e$ and length $\mathrm{l}(\gamma)=\nu-1$ satisfying the following conditions:
\begin{itemize}
\item[a)] the Littlewood-Richardson coefficient $c^\gamma_{\alpha \beta}$ is non-zero,
\item[b)] for every integer $2\leq i\leq \nu-1$, $\gamma_1> \gamma_i$,
\item[c)] $\sum_{i=1}^{\nu-1} (\gamma_i-2)^2=2(\sum_{j=1}^m \binom{\alpha_j}{2}+\sum_{k=1}^n \binom{\beta_k}{2})-4(e-\nu+1)$,
\item[d)] $\sum_{i=1}^{\nu-1} (\gamma_i-2)^3=6(\sum_{j=1}^m \binom{\alpha_j}{3}+\sum_{k=1}^n \binom{\beta_k}{3})+8(e-\nu+1)$.
\end{itemize}
\end{definition}

For example, if $\alpha=(3)$ and $\beta=(1,1,1)$, then it is easy to see that $\mathrm{P}(\alpha,\beta)=\{(4,1,1)\}$. Now we are ready to state our main result.

\begin{theorem}\label{main}
Let $\Gamma$ be a bipartite graph with colour classes $X$ and $Y$. Also assume that $\alpha=(\alpha_1, \dots, \alpha_m)$ and $\beta=(\beta_1, \dots, \beta_n)$ are weakly decreasing degree sequences of vertices of $X$ and $Y$, respectively. Set  $\nu:=m+n$ and $e:=\sum_{i=1}^m \alpha_i$. If the line graph $\mathrm{L}(\Gamma)$ of the graph $\Gamma$ is integral, then:
\begin{itemize}
\item[a)]  $\mathrm{P}(\alpha, \beta)$ is non-empty.
\item[b)]  There exists $\gamma \in \mathrm{P}(\alpha, \beta)$ such that the eigenvalues of $\mathrm{L}(\Gamma)$ are precisely $\gamma_1-2,\gamma_2-2,\dots,\gamma_{\nu-1}-2$, and $-2$ with multiplicity $e-\nu+1$.
\item[c)]  The diameter of $\mathrm{L}(\Gamma)$ is at most the maximum value of $\mathrm{k}(\gamma)$, where $\gamma$ runs over the set $\mathrm{P}(\alpha,\beta)$.
\end{itemize}
\end{theorem}

\begin{remark}
To compute the Littlewood-Richardson coefficient $c^\gamma_{\alpha\beta}$, the following are remarkable:
\begin{itemize}
\item[a)] There exists a polynomial time algorithm for deciding $c^\gamma_{\alpha\beta}$ is positive \cite{Bu2}.
\item[b)] There are some reduction formulae for Littlewood-Richardson coefficient (see \cite{Cho}, \cite{G} and \cite{Ki}).
\end{itemize}
\end{remark}

\begin{corollary}\label{integral}
Suppose $G$ is an integral line graph of a  bipartite graph.
Then the diameter of $G$ is at most $2\omega(G)$, where $\omega(G)$ is the clique number of $G$.
\end{corollary}

For a fixed integer $k\geq 3$, suppose $\Gamma$ is a $k$-regular graph. Letting $\lambda(\Gamma)$ be the second largest eigenvalue of $\Gamma$, it is a theorem of Alon and Boppana \cite{Alo} that $\lambda(\Gamma)\geq 2\sqrt{k-1}+\mathrm{O}(1)$, where $\mathrm{O}(1)$ goes to zero as $|\Gamma|\rightarrow \infty$. The graph $\Gamma$ is called a \textit{Ramanujan graph} if $|\lambda(\Gamma)|\leq 2 \sqrt{k-1}$. Lubotzky, Phillips and Sarnak \cite{LPS}, and independently Margulis \cite{Mar}, constructed the first examples of Ramanujan graphs; they are Cayley graphs of $\mathrm{PGL}_2(\mathbb{Z}/N\mathbb{Z})$ or  $\mathrm{PSL}_2(\mathbb{Z}/N\mathbb{Z})$ with $p+1$ explicit generators, for every prime $p$ and natural number $N$.
 Marcus, Spielman and Srivastava have proved the existence of $k$-regular bipartite Ramanujan graphs for arbitrary $k$ \cite{Mss}.

The \textit{bipartite complement} of a bipartite graph $\Gamma$ with two colour classes $X$ and $Y$ is the bipartite graph $\overline{\overline{\Gamma}}$ with the same colour classes having the edge between $X$ and $Y$ exactly where $\Gamma$ does not. The disjoint union of $n$ copies of the graph $\Gamma$ is denoted by $n\Gamma$. Also the disjoint union of two graphs $\Gamma_1$ and $\Gamma_2$ is denoted by $\Gamma_1 \cup \Gamma_2$.  Now we wish to present a restriction on the structure of a family of Ramanujan graphs.

 \begin{theorem}\label{ramanujan}
 Suppose $R$ is a non-complete Ramanujan graph. Then:
\begin{itemize}
\item[a)] If for some odd integer $k\geq 19$, $R$ is $k$-regular, then $R$ has an eigenvalue less than $-2$.
\item[b)] If $R$ is an integral line graph of a regular bipartite graph, then $R=\mathrm{L}(\Gamma)$, where $\Gamma$ is isomorphic to one of the following graphs:
\begin{enumerate}
\item[i)] $\Gamma\cong K_{s,s}$, for some positive integer $3\leq s\leq 10$.
\item[ii)] $\Gamma \cong \overline{\overline{(s+1)K_2}}$, for some positive integer $3\leq s\leq 8$.
\item[iii)] $\Gamma$ is isomorphic to one of the graphs $G_3-G_8$ described in \cite[fig 1]{Sc}.
\item[iv)] $\Gamma$ is isomorphic to one of the graphs $\overline{\overline{C_4\cup C_4\cup C_4}}$, $\overline{\overline{C_6\cup C_6}}$, $ \overline{\overline{C_4\cup C_4\cup C_6}}$, $\overline{\overline{C_4\cup C_4\cup C_4\cup C_4}}$, or $\overline{\overline{C_4\cup C_6\cup C_6}}$.
\item[v)] $\Gamma$ is isomorphic to one of the graphs $G_1,\overline{\overline{G_1}},G_2,\overline{\overline{G_2}},G_9,\overline{\overline{G_9}},G_{10},\overline{\overline{G_{10}}},G_{11},\,\\
 \overline{\overline{G_{11}}},\,G_{15}-G_{20},\,\overline{\overline{G_{15}}}-\overline{\overline{G_{20}}},\,G_{28}-G_{35},\,\overline{\overline{G_{34}}},\,\overline{\overline{G_{35}}},\,G_{37},\,\overline{\overline{G_{37}}},\,G_{39},\,\overline{\overline{G_{39}}},\,G_{41},\,\\
 G_{42}$ and $G_{43}$ described in \cite[Table 1]{Ko}.
 \end{enumerate}
 \end{itemize}
 \end{theorem}


\section{eigenvalues of line graphs}
 $\noindent$ In this section, we wish to prove our main results. A graph is called \textit{semi-regular bipartite}, with parameters $(n_1,n_2,r_1, r_2)$, if it is bipartite and vertices in the same colour class have the same degree ($n_1$ vertices of degree $r_1$ and $n_2$ vertices of degree $r_2$, where $n_1r_1=n_2r_2$).
   The identity matrix of rank $n$ is denoted by $I_n$. Also we use the notations $\mathrm{A}(\Gamma)$ and $p_{\Gamma}(x)$  for the adjacency matrix and the characteristic polynomial of a graph $\Gamma$, respectively. We begin with a useful observation on the eigenvalues of a  Hermitian matrix $C$ which is a consequence of Horn's inequalities.

 \begin{lemma}\label{fulton}\cite{Fulton}
Suppose $A$, $B$ and $C$ are Hermitian $n$ by $n$ matrices with $C=A+B$. Also assume that the eigenvalues of $A$ is denoted by $\alpha:=(\alpha_1, \alpha_2,\dots, \alpha_n)$, where $\alpha_1\geq \alpha_2\geq\dots\geq \alpha_n$. Similarly write $\beta$ and $\gamma$ for the eigenvalues  of $B$ and $C$, respectively. Then
  \begin{equation}
	\operatornamewithlimits{Max}_{i+j = n+k} \alpha_i + \beta_j \, \leq \,
\gamma_k \, \leq \, \operatornamewithlimits{Min}_{i+j = k+1} \alpha_i +
\beta_j,
\end{equation}
 for every positive integer $1\leq k\leq n$.
\end{lemma}

 \noindent \textbf{Proof of Theorem \ref{main}:} It is well-known that the least eigenvalue of the graph $\mathrm{L}(\Gamma)$ is equal to, or greater than $-2$. Note that the multiplicity of the eigenvalue $-2$ is equal to $e-\nu+1$ (see \cite[Theorem 2.2.4]{spect}). Thus we can assume that the eigenvalues of $\mathrm{L}(\Gamma)$ are precisely $\lambda_1\geq \lambda_2\geq \dots\geq \lambda_{\nu-1}$ and $-2$ with multiplicity $e-\nu+1$.

 Suppose $\Gamma_\alpha$ (resp. $\Gamma_\beta$) is a subgraph of $\mathrm{L}(\Gamma)$ whose vertex set is the vertex set of $\mathrm{L}(\Gamma)$, and two vertices $e_1$ and $e_2$ are adjacent in $\Gamma_\alpha$ (resp. $\Gamma_\beta$), if they have a common end in the colour class $X$ (resp. $Y$). Obviously, $\mathrm{A}(\mathrm{L}(\Gamma))=\mathrm{A}(\Gamma_\alpha)+\mathrm{A}(\Gamma_\beta)$. For every $1\leq i\leq \nu-1$, set $\gamma_i:=\lambda_i+2$. Applying Cayley–Hamilton theorem, we deduce that $\{\gamma_1,\gamma_2,\dots,\gamma_{\nu-1},0,\dots,0\}$, $\{\alpha_1,\alpha_2,\dots,\alpha_{m},0,\dots,0\}$ and $\{\beta_1,\beta_2,\dots,\beta_{n},0,\dots,0\}$ are the set of all eigenvalues of matrices $\mathrm{A}(\mathrm{L}(\Gamma))+2I_e$, $\mathrm{A}(\Gamma_\alpha)+I_e$ and $\mathrm{A}(\Gamma_\beta)+I_e$, respectively.

We now set $\gamma:=(\gamma_1, \gamma_2,\dots,\gamma_{\nu-1})$. Since $\mathrm{A}(\mathrm{L}(\Gamma))+2I_e=(\mathrm{A}(\Gamma_\alpha)+I_e)+(\mathrm{A}(\Gamma_\beta)+I_e)$, Theorems \ref{horn} and \ref{good} imply that the Littlewood-Richardson coefficient $c^\gamma_{\alpha\beta}$ is positive. Also as $\mathrm{L}(\Gamma)$ is connected, by the Perron-Frobenius theorem, $\gamma_1>\gamma_i$ for all $i\geq 2$. Since the sum of the $k$-th powers of the eigenvalues is just the number of closed walks of length $k$, it is easy to see that $\sum_{i=1}^{\nu-1} (\gamma_i-2)^2=2(\sum_{j=1}^m \binom{\alpha_j}{2}+\sum_{k=1}^n \binom{\beta_k}{2})-4(e-\nu+1)$ and $\sum_{i=1}^{\nu-1} (\gamma_i-2)^3=6(\sum_{j=1}^m \binom{\alpha_j}{3}+\sum_{k=1}^n \binom{\beta_k}{3})+8(e-\nu+1)$. Hence $\gamma \in \mathrm{P}(\alpha,\beta)$ and the eigenvalues of $\mathrm{L}(\Gamma)$ are precisely $\gamma_1-2,\gamma_2-2,\dots,\gamma_{\nu-1}-2$, and $-2$ with multiplicity $e-\nu+1$.

 Suppose $\mathrm{L}(\Gamma)$ has exactly $m$ distinct eigenvalues. It is well-known that the diameter $\rm{diam(L}(\Gamma))$ of the graph $\mathrm{L}(\Gamma)$ is bounded by this number, i.e, $\rm{diam(L}(\Gamma))\leq m-1$ (see \cite[Theorem 3.13]{Cv}). Thus as $\gamma \in \mathrm{P}(\alpha,\beta)$, we deduce that the diameter of $\mathrm{L}(\Gamma)$ is at most the maximum value of $\mathrm{k}(\lambda)$, where $\lambda$ runs over the set $\mathrm{P}(\alpha,\beta)$.
 It completes the proof of Theorem \ref{main}.\qed\\

 \noindent \textbf{Proof of Corollary \ref{integral}:}
By assumption, there exists a bipartite graph $\Gamma$ such that $G=\mathrm{L}(\Gamma)$.
 Suppose $\Delta$ is the maximum degree of the graph $\Gamma$.
 Obviously, $\Delta=\omega(G)$.
  Let $\alpha=(\alpha_1, \dots, \alpha_m)$ and $\beta=(\beta_1, \dots, \beta_n)$ be weakly decreasing degree sequences of vertices of colour classes of $\Gamma$.
   Then by Theorem \ref{main} (a),  $\mathrm{P}(\alpha, \beta)$ is non-empty.
    Let $\gamma\in \mathrm{P}(\alpha, \beta)$.
     Since $c^\gamma_{\alpha \beta}$ is positive, $\gamma_1\leq 2\Delta$.
     Hence using Theorem \ref{main} (c), the diameter of $G$ is at most $2\Delta=2\omega(G)$. \qed\\

    \noindent \textbf{Proof of Theorem \ref{ramanujan}:}
\textbf{a)} On the contrary, we assume that $R$ is a graph with least eigenvalue greater than or equal to $-2$. Since $k \geq 19$ is an odd integer, using \cite[Theorem 2. 5]{doob}, \cite{nn} and \cite[Theorem 4.1.5]{spect}, we deduce that $R$ is a line graph. By \cite[Proposition 1.1.5]{spect}, $R=\mathrm{L}(\Gamma)$, where $\Gamma$ is either regular or semi-regular bipartite. If $\Gamma$ is regular, then $k$ must be even which is impossible. Hence we can assume that $\Gamma$ is a semi-regular bipartite graph with parameters $(n_1,n_2,r_1,r_2)$ and colour classes $X$ and $Y$. Let $r_1\leq r_2$.
 Also assume that $R_X$ (resp. $R_Y$) is a subgraph of $R$ whose vertex set is the vertex set of $R$, and two vertices $e_1$ and $e_2$ are joined by an edge in $R_X$ (resp. $R_Y$), if they have a common end in the colour class $X$ (resp. $Y$).
 Obviously, $\mathrm{A}(R)=\mathrm{A}(R_X)+\mathrm{A}(R_Y)$.
  Let $\lambda(R)$ be the second largest eigenvalue of $R$. Then as $R$ is a non-complete Ramanujan graph, using Lemma \ref{fulton}, we deduce that $r_2-2\leq \lambda(R)\leq 2\sqrt{r_1+r_2-3} \leq 2\sqrt{2r_2-3}$.
   Thus $3\leq r_2 \leq 10$. Hence $k=r_1+r_2-2\leq 18$ which is a contradiction.\\
 \noindent \textbf{b)} Since $R$ is the line graph of a regular bipartite graph, there exists a bipartite $s$-regular graph $\Gamma$ with colour classes $X$ and $Y$ such that $R=\mathrm{L}(\Gamma)$.  Suppose $\lambda(R)$ is the second largest eigenvalue of the graph $R$.
  If $s+1\leq \lambda(R)$, then as $R$ is Ramanujan, $(s+1)^2\leq 8s-12$ which is a contradiction. Hence as $R$ is integral, $\lambda(R)=s-2,s-1$ or $s$. Let $n:=|X|=|Y|$. Then by Theorem \ref{main}, there exist non-negative integers $x$ and $y$ such that the spectrum of $R$ is
  \begin{equation}
\mathrm{Spec}(R)=\{-2^{(s-2)n+1}, (s-4)^{x}, (s-3)^{y}, (s-2)^{2n-2x-2y-2}, (s-1)^{y}, s^{x}, 2s-2\}.
\end{equation}
 Thus by \cite[Theorem 1.2.16]{spect},
$$\mathrm{Spec}(\Gamma)=\{-s, -2^{x}, -1^{y}, 0^{2n-2x-2y-2}, 1^{y}, 2^{x}, s\}.$$

 Suppose $\lambda(\Gamma)$ is the second largest eigenvalue of the graph $\Gamma$.   Now one of the following cases occurs:\\
 \textbf{Case 1.} $\lambda(\Gamma)=0$ (resp.$\lambda(\Gamma)=1$). Then as $R$ is Ramanujan, $(s-2)^2\leq 8s-12$ (resp. $(s-1)^2\leq 8s-12$). Hence $3\leq s\leq 10$ (resp. $3\leq s\leq 8$), and using \cite[Theorem 3]{Kol}, $\Gamma \cong K_{s,s}$ (resp. $\Gamma \cong \overline{\overline{(s+1)K_2}}$).\\
 \textbf{ Case 2.}  $\lambda(\Gamma)=2$. Since $R$ is Ramanujan, $s^2\leq 8s-12$. Hence $3\leq s\leq 6$. If $s=3$, then using \cite{Sc}, $\Gamma$ is isomorphic to one of the graphs $G_3-G_8$ described in \cite[fig 1]{Sc}. Thus we can assume that $4\leq s\leq 6$. Let $n\leq s+2$. Then as $\frac{p_{\Gamma}(x)}{x^2-s^2}=\frac{p_{\overline{\overline{\Gamma}}}(x)}{x^2-(n-s)^2}$ and $\lambda(\Gamma)=2$, we deduce that $\overline{\overline{\Gamma}}$ is a disjoint union of integral cycles. Hence  $\Gamma$ is isomorphic to one of the graphs $\overline{\overline{C_4\cup C_4\cup C_4}}$, $\overline{\overline{C_6\cup C_6}}$, $ \overline{\overline{C_4\cup C_4\cup C_6}}$, $\overline{\overline{C_4\cup C_4\cup C_4\cup C_4}}$ and $\overline{\overline{C_4\cup C_6\cup C_6}}$. Thus we can assume that $s+3\leq n$.
Then using Propositions 3.3, 3.4 and 3.5 of \cite{Ko}, we deduce that  $\Gamma$ is isomorphic to one of the graphs $G_1,\overline{\overline{G_1}},G_2,\overline{\overline{G_2}},G_9,\overline{\overline{G_9}},G_{10},\overline{\overline{G_{10}}},G_{11},\,
 \overline{\overline{G_{11}}},\,G_{15}-G_{20},\,\overline{\overline{G_{15}}}-\overline{\overline{G_{20}}},\,G_{28}-G_{35},\,\overline{\overline{G_{34}}},\,\overline{\overline{G_{35}}},\,G_{37},\,\overline{\overline{G_{37}}},\,G_{39},\,\overline{\overline{G_{39}}},\,G_{41},\,G_{42}$ and $G_{43}$ described in \cite[Table 1]{Ko}. This completes the proof.  \qed

\section*{Declaration of interests}
I declare that I have no known competing financial interests or personal relationships that could have appeared to influence the work reported in this paper.

\section*{Acknowledgements}
This research was supported in part
by a grant  from School of Mathematics, Institute for Research in Fundamental Sciences (IPM).


\end{document}